\documentclass[12pt]{amsart}

\usepackage{epsf}
\usepackage{graphicx}
\usepackage{latexsym}
\usepackage{amsfonts}
\usepackage{amsmath}
\usepackage{amssymb}
\usepackage{amscd}
\usepackage{amsthm}
\usepackage{bm}

\def\N{\mathbb{N}}

\def\C{\mathbb{C}}
\def\Z{\mathbb{Z}}

\def\C {\mathcal{C}}
\def\D {\mathcal{D}}
\def\T {\mathcal{T}}

\def\H {\mathbb{H}}
\def\s {\mathfrak{s}}

\newcommand{\abs}[1] {\left\lvert #1 \right\rvert}

\newcommand{\gen}[1] {\left\langle #1 \right\rangle}

\DeclareMathOperator{\coker}{coker} 
 \DeclareMathOperator{\Spin}{Spin}
 \DeclareMathOperator{\Char}{Char}
 
\DeclareMathOperator{\HF}{HF} \DeclareMathOperator{\HFK}{HFK}

\renewcommand{\tilde}{\widetilde}
\renewcommand{\bar}{\overline}
\renewcommand{\hat}{\widehat}

\newtheorem{thm}{Theorem}

\theoremstyle{definition}

\theoremstyle{remark}

\def \co{\colon\thinspace}

    \newcounter{myfootertablecounter}

    \newcommand\myfootnotemark{%
    \addtocounter{footnote}{1}%
    \footnotemark[\thefootnote]%
    }%

    \newcommand\myfootnotetext[1]{%
    \addtocounter{myfootertablecounter}{1}
    \footnotetext[\value{myfootertablecounter}]{#1}
    }

\begin{document}
\title
 [On Knots with Infinite Smooth Concordance Order]
 {On Knots with Infinite Smooth Concordance Order}

\author{Adam Simon Levine}
\address {Department of Mathematics, Columbia University \\ New York, NY 10027}
\email {alevine@math.columbia.edu}

\begin{abstract}
We use the Heegaard Floer obstructions defined by Grigsby, Ruberman,
and Strle to show that forty-six of the sixty-seven knots through
eleven crossings whose concordance orders were previously unknown
have infinite concordance order.
\end{abstract}

\maketitle

Let $K$ be an oriented knot in $S^3$. If $K$ bounds a smoothly
embedded disk in $D^4$, we say that $K$ is \emph{(smoothly) slice}.
Two knots $K, K'$ are said to be \emph{(smoothly) concordant} if $K
\# \bar{K'}$ is slice, where $\bar{K'}$ denotes the mirror of $K'$.
The set of concordance classes of knots forms a group $\C_1$ under
the connect sum operation with identity the unknot. The
\emph{concordance order} of a knot $K$ is the order of $K$ in
$\C_1$. The structure of the torsion in $\C_1$ is of considerable
interest; see, for instance, Livingston-Naik \cite{LN1, LN2} and
Jabuka-Naik \cite{JN}.

Let $Y_K = \Sigma_2(K)$ be the double branched cover of $K$, and let
$\tilde K$ be the inverse image of $K$ in $Y$. Grigsby, Ruberman,
and Strle \cite{GRS} defined numerical invariants $\D_n(K)$ and
$\T_n(K)$ ($n \in \N$) coming from the Heegaard Floer homology of
$Y_K$ and the knot Floer homology of $\tilde K$. They proved:
\begin{thm} \label{thm:GRS}
Let $K$ be a knot in $S^3$. Let $p$ be prime, and suppose that $p^m$
is the largest power of $p$ that divides $\det(K)$. If $K$ has
finite concordance order, then for each integer $0 \le e \le
\left\lfloor \frac{m+1}{2} \right\rfloor$, we have $\D_{p^e}(K) =
\T_{p^e}(K) = 0$.
\end{thm}
In practice, we are usually interested in $\D_p(K)$ and $\T_p(K)$,
where $p$ is 1 or a prime that divides $\det(K)$, so we restrict our
discussion to this case.

According to Livingston's database \emph{KnotInfo} \cite{CL}, the
smooth concordance orders of sixty-seven knots with up to eleven
crossings, listed in Table \ref{table:unknown}, were previously
unknown. We show here that forty-six of these knots, listed in
Tables \ref{table:alt} and \ref{table:nonalt}, have at least one
nonzero $\D_p$ invariant and hence have infinite concordance order.
For the remaining knots, all of the relevant $\D_p$ invariants
vanish, so the concordance orders of these knots remains unknown.
The $\T_p$ invariants for several of these knots can be obtained
using the author's computations of $\hat\HFK(Y_K, \tilde K)$
\cite{L}, but we do not obtain any new concordance information in
this manner.

\begin{table}
\[
\begin{array}{llllllll}
9_{30} & 9_{33} & 9_{44} & 10_{58} & 10_{60} & 10_{91} & 10_{102} &
10_{119} \\
10_{135} & 10_{158} & 10_{164} & 11a_{4} & 11a_5 & 11a_{8} &
11a_{11} & 11a_{24} \\
11a_{26} & 11a_{30} & 11a_{38} & 11a_{44} & 11a_{47} & 11a_{52} &
11a_{56} & 11a_{67} \\
11a_{72} & 11a_{76} & 11a_{80} & 11a_{88} & 11a_{98} & 11a_{104} &
11a_{109} & 11a_{112} \\
11a_{126} & 11a_{135} & 11a_{160} & 11a_{167} & 11a_{168} &
11a_{170} & 11a_{187} & 11a_{189} \\
11a_{233} & 11a_{249} & 11a_{257} & 11a_{265} & 11a_{270} &
11a_{272} & 11a_{287} & 11a_{288} \\
11a_{289} & 11a_{300} & 11a_{303} & 11a_{315} & 11a_{350} & 11n_{12}
& 11n_{34} & 11n_{45} \\
11n_{48} & 11n_{53} & 11n_{55} & 11n_{85} & 11n_{100} & 11n_{110} &
11n_{114} & 11n_{130} \\
11n_{145} & 11n_{157} & 11n_{165}
\end{array}
\]
\caption{Knots through eleven crossings with unknown concordance
order. \label{table:unknown}}
\end{table}

\begin{table}
\begin{tabular}{|c|c|c|}
  \hline
  Knot $K$ & $\det(K)$ & Nonzero GRS invariants \\
  \hline
  $9_{30}$ & 53 & $\D_{53}=4$ \\
  $9_{33}$ & 61 & $\D_{61}=4$  \\
  $10_{58}$ & 65 & $\D_{13}=4$ \\
  $10_{60}$ & 85 & $\D_{17}=4$ \\
  $10_{102}$ & 73 & $\D_{73}=-12$ \\
  $10_{119}$ & 101 & $\D_{101}=-16$\\
  $11a_{4}$ & 97 & $\D_{97}=-24$ \\
  $11a_{8}$ & 117 & $\D_{13}=-4$ \\
  $11a_{11}$ & 113 & $\D_{113}=12$ \\
  $11a_{24}$ & 157 & $\D_{157}=12$ \\
  $11a_{26}$ & 157 & $\D_{157}=12$ \\
  $11a_{30}$ & 149 & $\D_{149}=12$ \\
  $11a_{52}$ & 137 & $\D_{137}=16$ \\
  $11a_{56}$ & 109 & $\D_{109}=-8$ \\
  $11a_{67}$ & 125 & $\D_{25}=-4$ \\
  $11a_{76}$ & 145 & $\D_{29}=-4$ \\
  $11a_{80}$ & 137 & $\D_{137}=-12$ \\
  $11a_{88}$ & 101 & $\D_{101}=-8$ \\
  $11a_{126}$ & 145 & $\D_5 = 4, \D_{29}=4$ \\
  $11a_{160}$ & 145 & $\D_{29}=-4$ \\
  $11a_{167}$ & 113 & $\D_{113}=12$ \\
  $11a_{170}$ & 185 & $\D_{37}=-4$ \\
  $11a_{189}$ & 149 & $\D_{149}=-12$ \\
  $11a_{233}$ & 173 & $\D_{101}=16$ \\
  $11a_{249}$ & 117 & $\D_{13}=-4$ \\
  $11a_{257}$ & 97 & $\D_{97}=-8$ \\
  $11a_{265}$ & 109 & $\D_{109}=24$ \\
  $11a_{270}$ & 137 & $\D_{137}=12$ \\
  $11a_{272}$ & 149 & $\D_{149}=12$ \\
  $11a_{287}$ & 181 & $\D_{181}=-12$ \\
  $11a_{288}$ & 205 & $\D_{5}=4, \D_{41}=4$ \\
  $11a_{289}$ & 145 & $\D_{29}=4$ \\
  $11a_{300}$ & 153 & $\D_{17}=-4$ \\
  $11a_{303}$ & 149 & $\D_{149}=36$ \\
  $11a_{315}$ & 157 & $\D_{157}=12$ \\
  $11a_{350}$ & 185 & $\D_{5}=4, \D_{37}=4$ \\
  \hline
\end{tabular}

\caption{Alternating knots with non-vanishing Grigsby-Ruberman-Strle
$\D_p$ invariants. \label{table:alt}}
\end{table}

\begin{table}
\begin{tabular}{|c|c|c|}
\hline
  Knot $K$ & $\det(K)$ & Nonzero GRS invariants \\
  \hline
  $9_{44}$ & 17 &  $\D_{17}=4$ \\
  $10_{135}$ & 135 & $\D_{37}=4$ \\
  $11n_{12}$ & 13 & $\D_{13}=-8$ \myfootnotemark \\
  $11n_{48}$ & 29 & $\D_{29}=-8$ \\
  $11n_{53}$ & 37 & $\D_{37}=-8$ \\
  $11n_{55}$ & 61 & $\D_{61}=12$ \\
  $11n_{110}$ & 41 & $\D_{41}=-12$ \\
  $11n_{114}$ & 53 & $\D_{53}=-4$ \\
  $11n_{130}$ & 53 & $\D_{53}=12$ \\
  $11n_{165}$ & 85 & $\D_{17}=-4$ \\
\hline
\end{tabular}
\caption{Non-alternating knots with non-vanishing
Grigsby-Ruberman-Strle $\D_p$ invariants. \label{table:nonalt}}
\end{table}

For the remainder of this paper, we describe the techniques used to
compute the $\D_p$ and $\T_p$ invariants for the knots considered
here.

Let us briefly recall the definition of these invariants in the case
where $H^2(Y_K;\Z)$ is cyclic. (For the general case, see
\cite[Definition 4.1]{GRS}.) Let $\s_0\in \Spin^c(Y_K)$ be the
so-called \emph{canonical spin$^c$ structure} on $Y_K$, uniquely
characterized by the property that $c_1(\s_0)=0$. Recall that
$\Spin^c(Y_K)$ is an affine space for $H^2(Y_K;\Z)$, so we may
identify $\Spin^c(Y_K)$ with $H^2(Y_K;\Z)$ via the identification
$\s \mapsto c_1(\s)$. Let $G_m$ be the unique order-$p$ subgroup of
$H^2(Y_K;\Z)$. The invariants $\D_p(K)$ and $\T_p(K)$ are then
defined as
\[
\begin{split}
\D_p(K) &= \sum_{\s \in \s_0 + G_p} d(Y_K,\s) \\
\T_p(K) &= \sum_{\s \in \s_0 + G_p} \tau(Y_K, \tilde K, \s).
\end{split}
\]
Here $d(Y_K,\s)$ is the correction term for $\HF^+(Y_K,\s)$, and
$\tau(Y_K,\tilde K,\s)$ is the $\tau$-invariant for $\hat\HFK(Y_K,
\tilde K, \s)$. (See Ozsv\'ath-Szab\'o \cite{OSzAbs, OSzKnot} for
the definitions of $d$ and $\tau$.)

In many cases, the results of Ozsv\'ath and Szab\'o
\cite{OSzPlumbed, OSzDouble, OSzUnknot} may be used to compute the
correction terms $d(Y,\s)$ combinatorially. Given a projection of
$K$, let $G$ be its Goeritz matrix (defined in \cite[section
3]{OSzDouble}). Let $\abs{G}$ denote the rank of $G$. The double
cover $Y_K$ bounds a 4-manifold $X_G$ whose intersection form on
$H_2$, $Q=Q_{X_G}$, is given by $G$ (with respect to a basis of
spheres). Let $\Char(G) \subset H^2(X_G;\Z)$ denote the set of
characteristic vectors for $Q$, i.e., vectors $\alpha \in
H^2(X_G;\Z)$ such that $\gen{\alpha,v} \equiv Q(v,v) \pmod 2$ for
every $v \in H_2(X_G;\Z)$. The restriction map $i^*\co H^2(X_G) \to
H^2(Y_K)$ partitions $\Char(G)$ into equivalence classes
$\Char(G,\s)$ corresponding to the spin$^c$ structures on $Y_K$.
Given certain hypotheses on $G$, including that $G$ is
negative-definite, Ozsv\'ath and Szab\'o \cite[Corollary
1.5]{OSzPlumbed} proved that the correction terms for $\HFK^+(Y_K)$
are given by the formula
\begin{equation} \label{eqn:dformula}
d(Y_K,\s) = \max_{\alpha \in \Char(G,\s)} \frac{\alpha^2 + \abs{G}}
{4}.
\end{equation}
Ozsv\'ath and Szab\'o provide an algorithm for finding the vectors
in each equivalence class that realize this maximum. Moreover, since
$\H^2(Y_K;\Z) \cong \coker(G)$, we may easily identify the group
structure on $\Spin^c(Y_K)$ (specifically, which spin$^c$ structures
are in the special subgroup $G_p$) using the Smith normal form for
$G$.

\myfootnotetext{That $K=11n_{12}$ has infinite concordance order
  also follows from the simpler fact that $\tau(S^3,K)=1$, as was computed by
  Baldwin and Gillam \cite{BG}.}

As shown in \cite{OSzDouble}, Equation \ref{eqn:dformula} holds
whenever $G$ is computed from an alternating projection. More
generally, if $K$ admits a projection that is alternating except in
a region that consists of left-handed twists, Ozsv\'ath and Szab\'o
\cite{OSzUnknot} show how to use Kirby calculus on $X_G$ to obtain a
matrix $\tilde G$ for $Q$ that satisfies the correct hypotheses.
(See also Jabuka-Naik \cite{JN} for a concise explanation.) All of
the non-alternating knots in Table \ref{table:unknown} satisfy this
hypothesis, so we may compute the $\D_p$ invariants as described
above.

Finally, to compute the $\T_p$ invariants of a knot, one must
compute the integers $\tau(Y_K,\tilde K, \s)$ associated to the
spectral sequence from $\hat\HFK(Y_K, \tilde K, \s)$ to
$\hat\HF(Y_K,\s)$. When $\hat\HFK(Y_K, \tilde K, \s)$ and
$\hat\HF(Y_K,\s)$ are sufficiently simple, one can sometimes
determine $\tau$ without knowing all the differentials in the
spectral sequence. For instance, if $\hat\HF(Y_K,\s)$ has rank 1 and
$\hat\HFK(Y_K, \tilde K, \s)$ is supported on a single diagonal,
$\tau(Y_K, \tilde K, \s)$ is equal to the Alexander grading of the
nonzero group in Maslov grading $d(Y_K,\s)$. The author \cite{L} has
shown how to compute $\HFK(Y_K, \tilde K)$ (with coefficients in
$\Z/2$) for any knot $K$ using grid diagrams and has computed the
values of $\tau$ for several of the non-alternating knots considered
here ($9_{44}$, $10_{135}$, $10_{158}$, $10_{164}$, $11n_{100}$, and
$11n_{145}$). However, the $\T_p$ invariants all vanish in these
cases, so we do not obtain any new concordance information.

\end{document}